\documentclass[10pt]{amsart}
\usepackage{amssymb}
\usepackage[centertags]{amsmath}
\usepackage{amssymb}
\usepackage{amscd}
\usepackage{amsthm}
\usepackage{amsxtra}
\newtheorem{thm}{Theorem}[section]
\newtheorem{prop}[thm]{Proposition}
\newtheorem{lemma}[thm]{Lemma}

\begin{document}
\title{Hochschild Cohomology of $\rm{II_1}$ Factors with cartan masas}

\author{Jan M. Cameron}
\address{Department of Mathematics\\ Texas A \& M University\\
College Station Texas, 77843-3368} \email{jcameron@math.tamu.edu}


\begin{abstract} In this paper we prove that for a type $\rm{II_1}$ factor $N$ with a
Cartan maximal abelian subalgebra (masa), the Hochschild
cohomology groups $H^n(N,N)$=0, for all $n \geq 1$.  This
generalizes the result of Sinclair and Smith, who proved this for
all $N$ having separable predual.
\end{abstract}

\maketitle

\section{Introduction}
     The study of Hochschild cohomology for von Neumann
algebras can be traced back to a well-known theorem, due
separately to Kadison \cite{Kadison derivations} and Sakai
\cite{Sakai derivations}, which states that every derivation
$\delta:M \rightarrow M$ on a von Neumann algebra $M$ is inner,
that is, there exists an element $a \in M$ such that
$\delta(x)=xa-ax,$ for all $x \in M$.  This corresponds to the
vanishing of the first continuous Hochschild cohomology group,
$H^1(N,N)$.  It is natural to conjecture that the higher
cohomology groups $H^n(N,N)$ are also trivial. This program was
taken up in the seventies, in a series of papers by Johnson,
Kadison, and Ringrose (\cite{JKR cohomology III, KR cohomology I,
KR cohomology II}), who affirmed the conjecture for type I
algebras and hyperfinite algebras.  In the mid 1980's, a parallel
theory of completely bounded cohomology was initiated. The groups
$H^n_{cb}(N,N)$ are computed under the additional assumption that
all cocycles and coboundaries, usually assumed to be norm
continuous, are completely bounded. Christensen and Sinclair
showed that $H^n_{cb}(N,N)=0$ for all von Neumann algebras $N$(see
\cite{cohomology book} for details of the proof). It was proved in
\cite{completely bdd} that continuous and completely bounded
cohomology coincide for all von Neumann algebras stable under
tensoring with the hyperfinite $\rm{II_1}$ factor $R$. Thus, it is
known that the continuous Hochschild cohomology vanishes for all
von Neumann algebras of type $\rm{I},$ $\rm{II_{\infty}}$, and
$\rm{III}$, and for all type $\rm{II_1}$ algebras stable under
tensoring with $R$. This leaves open the case of the general type
$\rm{II_1}$ von Neumann algebra. By direct integral techniques, it
is enough in the separable case to compute the groups $H^n(N,N)$
when $N$ is a factor.

     As in the Kadison-Sakai theorem, the vanishing of Hochschild
cohomology groups gives structural information about certain
bounded $n$-linear maps on the von Neumann algebra, in particular,
that they can be formed from the $(n-1)$-linear maps.  The
vanishing of certain other higher cohomology groups also yields
some nice perturbation results for von Neumann algebras (see
\cite[Chapter 7]{cohomology book}). The Hochschild cohomology
groups have been shown to be trivial for a few large classes of
type $\rm{II_1}$ factors, including those with property $\Gamma$,
\cite{property gamma}, and those which contain a Cartan maximal
abelian subalgebra (masa) and have separable predual, \cite{cartan
sep}. The purpose of this paper is to extend this last result to
include arbitrary $\rm{II_1}$ factors with a Cartan masa.  The
techniques in \cite{cartan sep} depended heavily on separability,
and so could not be modified to encompass nonseparable algebras.
The techniques in the present note originated in \cite{property
gamma} and \cite{strongly singular}, and can be viewed as part of
a general strategy to relate properties of nonseparable type
$\rm{II_1}$ factors to their separable subalgebras.

    The author thanks Roger Smith for his patient guidance
throughout the completion of this work.

\section{Notation and Background}\label{section:preliminaries}

\subsection{Finite von Neumann algebras}\label{subsection:finite}

The setting of this paper is a finite von Neumann algebra $N$ with
a faithful, normal trace $\tau$. We will call such an algebra
\emph{separable} if one of the following equivalent conditions
holds \cite{takesaki}:
\begin{itemize}
\item[\rm(i)] There exists a countable set of projections in $N$
generating a weakly dense subalgebra of $N$.

\item[\rm(ii)] The Hilbert space $L^2(N,\tau)$ (on which $N$ is
faithfully represented) is separable in its norm.

\item[\rm(iii)] The predual of $N$ is a separable Banach space.

\end{itemize}

The weakly closed unit ball of $N$ is then separable in this
topology, being a compact metric space.

     If $A\subseteq M\subseteq N$  is an inclusion of von Neumann
subalgebras, we denote the \emph{normalizer of A in M} by
\[\mathcal{N}(A,M) = \{u\in \mathcal{U}(M): uAu^*=A\},\]
where $\mathcal{U}(M)$  denotes the unitary group of $M$. We will
write $\mathcal{N}(A)$ for $\mathcal{N}(A,N)$, and refer to this
set as the \emph{normalizer} of $A$.  Dixmier, \cite{Dixmier},
classified the masas $A$ of a von Neumann algebra $M$ as singular,
Cartan, or semi-regular according to whether
$\mathcal{N}(A,M)''$=$A$, $\mathcal{N}(A,M)''=M$, or
$\mathcal{N}(A,M)''$ is a proper subfactor of $M$. See the
forthcoming notes of Sinclair and Smith, \cite{masabook}, for a
comprehensive exposition of the theory of masas in von Neumann
algebras.

    The following facts about finite von Neumann algebras are
standard, but are critical to what follows, so we include some
discussion for the reader's convenience.  Recall that if $Q$ is a
von Neumann subalgebra of a finite von Neumann algebra $N$ with
trace $\tau$, then there exists a bounded conditional expectation
$\mathbb{E}_Q:N \rightarrow Q$ which is characterized by the
following two properties:
\begin{itemize}
\item[\rm{(i)}] For all $x\in N$ and $q_1$, $q_2 \in Q$,
$\mathbb{E}_Q(q_1xq_2)=q_1\mathbb{E}_Q(x)q_2$ ($Q$-bimodularity)

\item[\rm{(ii)}] For all $x\in N$, $\tau(\mathbb{E}_Q(x))=\tau(x)$
($\mathbb{E}_Q$ is trace-preserving).
\end{itemize}
We emphasize that $\mathbb{E}_Q$ is unique among all maps from $N$
into $Q$ with these two properties, and not just those assumed to
be continuous and linear.  This observation leads to the following
result, which is essentially in \cite{subalgebras}, and can be
found in \cite{masabook}.

\begin{prop}\label{closure}
Let $N$ be a finite von Neumann algebra, and let $Q$ be a von
Neumann subalgebra of $N,$ with faithful, normal trace $\tau$. For
any $x\in N,$ $\mathbb{E}_{Q'\cap N}(x)$ is the unique element of
minimal $\|\cdot\|_2$-norm in the weak closure of
\[K_Q(x)= {\rm{conv}} \{uxu^*:u \in \mathcal{U}(Q)\}.\]
\end{prop}

It is important to consider whether the weak closure of the set
$K_Q(x)$ is taken in $N$ or $L^2(N)$.  The embedding of $N$ into
$L^2(N)$ is continuous, when both spaces are given their
respective weak topologies.  Thus, the (compact) weak closure of
any ball in $N$ is weakly closed (hence also $\|\cdot\|_2$-closed,
by convexity) in $L^2(N)$.  Conversely, the preimage of a weakly
closed (equivalently, $\|\cdot\|_2$-closed) ball in $L^2(N)$ is
weakly closed in $N$.  Thus, the weak closure in $N$, the weak
closure in $L^2(N)$, and the $\|\cdot\|_2$ closure in $L^2(N)$ of
$K_Q(x)$ all coincide.  We will use the common notation $K_Q^w(x)$
for all three closures.

\begin{lemma}\label{lemma:slippery}
Let N be a finite von Neumann algebra with a normal, faithful
trace $\tau$.  Let M be a separable subalgebra of N.  Suppose
there exists a subset $S\subseteq N$ of unitaries such that
$S''=N$. Then there exists a countable subset F of S such that
$F''\supseteq M$.

\end{lemma}
\begin{proof}  We first claim that $C^*(S)$ is $\|\cdot\|_2$-dense in
$L^2(N).$  Let $x\in N$.  Since $C^*(S)$ is strongly dense in $N$,
by the Kaplansky density theorem, there exists a net
$\{x_{\alpha}\}$ in $C^*(S)$ such that $\|x_{\alpha}\| \leq \|x\|$
and $x_{\alpha}$ converges to $x$ $*$-strongly.  Then also
$(x_{\alpha}-x)^*(x_{\alpha}-x)$ converges to 0 weakly.  Moreover,
$\|(x_{\alpha}-x)^*(x_{\alpha}-x)\|$ is uniformly bounded and
$\tau$ is a normal state (hence weakly continuous on bounded
subsets of $N$ by \cite[Theorem 7.1.12]{KR}); then
$\|x_{\alpha}-x\|_2^2 = \tau((x_{\alpha}-x)^*(x_{\alpha}-x))$
converges to zero.  Thus, $x \in \overline{C^*(S)}^{\|\cdot\|_2},$
and the claim follows.

      Note that $L^2(M)$ is a separable Hilbert subspace of
$L^2(N)$.  Let $\{\xi_n\}_{n=1}^\infty$ be a dense subset.  For
each $n$, there is a sequence $\{s_{nk}\}_{k=1}^\infty$ in
$C^*(S)$ such that $s_{nk}$ converges in $\|\cdot\|_2$-norm to
$\xi_n$. The operators $s_{nk}$ lie in the norm closure of
${\rm{Alg}}(F)$, for some countable subset $F$ of $S$.  It follows
that $L^2(M) \subseteq L^2({\rm{Alg}}(F)) \subseteq L^2(F'')$.  We
claim that this implies $M\subseteq F''$.  Let
$\mathbb{E}_M:N\rightarrow M$ and $\mathbb{E}_{F''}:N\rightarrow
F''$ denote the respective trace-preserving conditional
expectations, obtained by restricting the Hilbert space
projections $e_M:L^2(N)\rightarrow L^2(M)$ and $e_{F''}:L^2(N)
\rightarrow L^2(F'')$ to $N.$  For any $x\in M$, we have
\[x=e_{F''}(x)+(1-e_{F''})(x)=e_{F''}(x),\] since $L^2(M)
\subseteq L^2(F'')$.  But then $x=\mathbb{E}_{F''}(x)$, so
$M\subseteq F''$.
\end{proof}
\subsection{Cohomology}\label{subsection:cohomology}
In this section we collect the basic definitions and results of
Hochschild cohomology that will be used in the next section.  The
reader may wish to consult \cite{cohomology book} for a more
detailed exposition. Let $M$ be a von Neumann algebra and let $X$
be a Banach $M$-bimodule (In the next section, we will restrict
attention to the case $X=M$).  Let $\mathcal{L}^n(M,X)$ denote the
vector space of $n$-linear bounded maps $\phi:M^n \rightarrow X$.
Define the coboundary map $\partial:\mathcal{L}^n(M,X)\rightarrow
\mathcal{L}^{n+1}(M,X)$ by
\begin{align*}
\partial \phi(x_1,\ldots,x_{n+1}) &= x_1 \phi (x_2,\ldots,x_{n+1})\\
                               &+  \sum_{j=1}^n (-1)^j \phi
                              (x_1,\ldots,x_{j-1},x_jx_{j+1},\ldots,x_{n+1})\\
                                &+ (-1)^{n+1} \phi(x_1,\ldots,x_n)x_{n+1}.
\end{align*}
An algebraic computation shows that $\partial^2=0$.  We thus
obtain the Hochschild complex and define the $n$th continuous
Hochschild cohomology group to be \[H^n(M,X)=\frac{ker
\partial}{Im \partial}.\]
 The maps $\phi \in ker\partial$ are called $(n+1)$-cocycles, and the maps in $Im
\partial$ are called $n$-coboundaries.  It is easy to check that the $2$-cocycles are
precisely the bounded derivations from $M$ into $X$.

     Various extension and averaging arguments are of central importance in computing
cohomology groups.  The most basic extension argument states that
normal cohomology and continuous cohomology are equal when the
space $X$ is assumed to be a dual normal $M$-module (see
\cite[Theorem 3.3]{cohomology survey}).  Thus, we will assume in
what follows that all $n$-cocycles are separately normal in each
variable. We denote the vector space of bounded, separately normal
$n$-linear maps from $M$ to $X$ by $\mathcal{L}^n_w(M,X)$. The
general strategy of the averaging arguments is to replace a
continuous $n$-cocycle $\phi$ with a modified cocyle $\phi +
\partial\psi$ which has more desirable continuity and modularity
properties.  The coboundary $\partial \psi$ is obtained by an
averaging process over a suitable group of unitaries in the
underlying von Neumann algebra.  In the present work, it will be
essential that modifications to a given cocycle are made while
preserving certain norm estimates on the original cocycle.  The
averaging result we need is the following, which can be found in
\cite{cohomology book}:

\begin{lemma} \label{normal R modular}  Let $R$ be a hyperfinite
von Neumann subalgebra of a von Neumann algebra $M$.  Then there
is a bounded linear map $L_n:\mathcal{L}^n_w(M,M) \rightarrow
\mathcal{L}^{n-1}_w(M,M)$ such that $\phi+ \partial L_n \phi$ is a
separately normal $R$-module map for any $n$-cocycle $\phi$.
Moreover, $\|L_n\| \leq \ell (n)$, a constant depending only on
$n$.
\end{lemma}
We will also need two extension results, the first of which is
essentially Lemma 3.3.3 of \cite{cohomology book}.  The second is
a modification of Lemma 3.3.4 in \cite{cohomology book}, as we
will not need the full generality of that result.
\begin{lemma} \label{extension}
Let $A$ be a $C^*$-algebra on a Hilbert space $\mathcal{H}$ with
weak closure $\overline{A}$ and let $X$ be a dual Banach space. If
$\phi:A^n \rightarrow X$ is bounded, $n$-linear, and separately
normal in each variable then $\phi$ extends uniquely, without
changing the norm, to a bounded, separately normal, $n$-linear map
$\overline{\phi}:(\overline{A})^n \rightarrow X$.
\end{lemma}

\begin{lemma} \label{extension2}
Let $A$ be a $C^*$-algebra, and denote its weak closure by
$\overline{A}$.  Then for each $n$ there exists a linear map
$V_n:\mathcal{L}^n_w(A,X) \rightarrow
\mathcal{L}^n_w(\overline{A},X)$ such that when $\phi$ is as in
Lemma \ref{extension}, then $V_n \phi = \overline{\phi}$.
Moreover, $\|V_n\| \leq 1$ and $\partial V_n= V_{n+1} \partial$,
for all $n \geq 1$.

\end{lemma}
\begin{proof}  For $\phi$ in $\mathcal{L}^n_w(A,X)$, we define
$V_n \phi = \overline{\phi}$.  This is well-defined and linear by
the uniqueness in Lemma \ref{extension}.  Moreover, since the
extension in Lemma \ref{extension} is norm-preserving, we have
$\|V_n\| \leq 1$.  Now for any $x_1,x_2,\ldots x_{n+1} \in A$, we
have
\[\partial V_n \phi(x_1,\ldots x_{n+1})=\partial
\phi(x_1,\ldots x_{n+1})=V_{n+1} \partial \phi
(x_1,\ldots,x_{n+1}).\] Both of the maps $\partial V_n \phi$ and
$V_{n+1}
\partial \phi$ are separately normal, so are uniquely defined on
$(\overline{A})^n$. Thus, we will have \[\partial V_n \phi
(x_1,\ldots x_{n+1})=V_{n+1} \partial \phi (x_1,\ldots
,x_{n+1}),\] for all $x_1,\ldots,x_{n+1}$ in $\overline{A}$.  This
completes the proof.
\end{proof}
\section{Main Results}\label{section:results}
    A corollary of the following result is that the study of
Cartan masas can, in many instances, be reduced to the separable
case.  The techniques of the proof come from \cite[Theorem
2.5]{strongly singular}, in which a similar result is proved for
singular masas in $\rm{II_1}$ factors.
\begin{prop}\label{prop:cartansep}
Let $N$ be a $\rm{II_1}$ factor with Cartan masa $A$, and let
$M_0$ be a separable von Neumann subalgebra of $N$.  Let
$\phi:N^n\rightarrow N$ be a separately normal $n$-cocycle. Then
there exists a separable subfactor $M$ such that $M_0\subseteq
M\subseteq N$, $M\cap A$ is a Cartan masa in $M$, and $\phi$ maps
$M^n$ into $M$.

\end{prop}
\begin{proof} For a von Neumann algebra $Q$, and $x \in Q,$
denote by $K_Q^n(x)$ and $K_Q^w(x)$, respectively, the operator
and $\|\cdot\|_2$-norm closures of the set $K_Q(x)$.  Recall that
when $Q$ is a von Neumann subalgebra of $M$, for any $x \in M$,
$\mathbb{E}_{Q' \cap M}(x)$ picks out the element of minimal
$\|\cdot\|_2$-norm in $K_Q^w(x)$. We will construct, inductively,
a sequence of separable von Neumann algebras
      \[M_0\subseteq M_1\subseteq M_2\subseteq\ldots\subseteq N\]
and abelian subalgebras $B_k\subseteq M_k$ so that
$M=(\bigcup_{k=1}^\infty M_k)''$ has the required properties.  The
von Neumann algebra $B=(\bigcup_{n=0}^{\infty}B_n)''$ will be a
masa in $M$, and thus equal to $M \cap A$.  The inductive
hypothesis is as follows:

For a fixed sequence $\{y_{k,r}\}_{r=1}^\infty$,
$\|\cdot\|_2$-norm dense in the $\|\cdot\|_2$-closed unit ball of
the separable von Neumann algebra $M_k$,
\begin{itemize}
\item[\rm{(i)}] $\mathbb{E}_A(y_{k,r})\in B_k\cap
K_{B_{k+1}}^w(y_{k,r})$,
 for $r\geq 1,$ where $B_k = M_k \cap A$;

 \item[\rm{(ii)}]$K_{M_{k+1}}^n(y_{k,r}) \cap \mathbb{C}1$ is nonempty for all
 $r \geq 1$;

 \item[\rm{(iii)}] For each k, there is a countable set of unitaries
 $\mathcal{U}_{k+1} \subseteq N(A) \cap  M_{k+1}$ such that
 $\mathcal{U}_{k+1}'' \supseteq M_k$;

 \item[\rm{(iv)}]  The cocycle $\phi$ maps $(M_k)^n$ into $M_{k+1}$.

\end{itemize}

We first prove that this sequence of algebras gives the desired
result.  The von Neumann algebra $M=(\bigcup_{n=0}^\infty M_n)''$
will be a separable subalgebra of $N$, by our construction.  We
show M is a factor.

     For any $x \in M,$ $K_M^w(x) \cap \mathbb{C}1$
is nonempty, by the following approximation argument.  First
suppose that $x \in M_k,$ for some $k \geq 1$ and $\|x\| \leq 1$.
Let $\varepsilon > 0$ be given.  Choose an element $y_{k,r} \in
M_k$ as above with $\|y_{k,r}-x\|_2 < \varepsilon$.  By condition
(ii) we can choose an element
\[a_{k,r}=\sum_{i=1}^m \lambda_i u_i y_{k,r} u_i^* \in
K_{M_{k+1}}(y_{k,r})\] whose $\|\cdot\|_2$-norm distance to
$\mathbb{C}1,$ which we denote
${\rm{dist}}_2(a_{k,r},\mathbb{C}1),$ is less than $\varepsilon$.
Then $a=\sum_{i=1}^m \lambda_i u_i x u_i^*$ is an element of
$K_{M_{k+1}}(x)$ with ${\rm{dist}}_2(a,\mathbb{C}1)<2
\varepsilon$. It follows that $K_{M_{k+1}}^w(x) \cap \mathbb{C}1$
is nonempty.

    Now let $x \in M,$ $\|x\| \leq 1$, and fix $\varepsilon >0$.  By the Kaplansky density
theorem, there exists a $k \geq 1,$ and an element $x_{\alpha} \in
M_k$ of norm at most 1 such that $\|x_{\alpha}-x\|_2 <
\varepsilon$. Then by what we did above,
$K_{M_{k+1}}^w(x_{\alpha}) \cap \mathbb{C}1$ is nonempty. It
follows, by similar argument to the one above, that there exists
an element $a \in K_M(x)$ with ${\rm{dist}}_2(a, \mathbb{C}1) <
2\varepsilon.$  Thus, $K_M^w(x)$ has nonempty intersection with
$\mathbb{C}1$.  By scaling, this result is true for $x \in M$ of
arbitrary norm.  To see that $M$ is a factor, note that if $x$ is
central in $M$ then $K_M^w(x)=\{x\}$, and since this set meets
$\mathbb{C}1$, $x$ must be a multiple of the identity.

     We now prove that $B=(\bigcup_{n=0}^{\infty}B_n)''$ is a masa in $M$. First,
condition (i) implies that $\mathbb{E}_A(x) \in B \cap K_B^w(x)
\subseteq B \cap K_A^w(x),$ for all $x \in M$.  This follows from
an approximation argument similar to the one above, in which we
prove the claim first for all $x \in M_k$ with $\|x\|\leq1 $ and
then extend to all of $M$. Since $\mathbb{E}_A(x)$ is the element
of minimal $\|\cdot\|_2$-norm in $K_A^w(x),$ it also has this
property in $K_B^w(x)$.  But then $\mathbb{E}_A(x) =
\mathbb{E}_{B' \cap M}(x).$  Since for all $x \in M$,
$\mathbb{E}_A(x)=\mathbb{E}_B(x)$, one has \[x=\mathbb{E}_{B' \cap
M}(x) = \mathbb{E}_A(x)=\mathbb{E}_B(x),\] for all $x \in B' \cap
M$.  Thus, $B' \cap M \subseteq B$.  Since $B$ is abelian, the
opposite inclusion also holds.  Thus $B$ is a masa in $M,$ and
$B=M \cap A$.  We show $B$ is Cartan.  By condition (iii), we will
have $M = (\bigcup_{k=0}^\infty \mathcal{U}_{k+1}'')''$.  We claim
that this last set is precisely $\mathcal{N}(B,M)''$, and hence
that $B$ is Cartan in $M$. Fix $k \geq 0$ and let $u \in
\mathcal{U}_{k+1}.$ Then since $u \in \mathcal{N}(A) \cap
M_{k+1}$, for any $j\leq k+1$ we have
\[uB_ju^*=u(A \cap M_j)u^* \subseteq A \cap M_{k+1} =
B_{k+1} \subseteq B.\] If $j>k+1$, since $u \in M_{k+1} \subseteq
M_{j},$ we have
\[uB_ju^*=u(A \cap M_j)u^* = A \cap M_j = B_j \subseteq B.\]

Then $u( \bigcup_{j=0}^\infty B_j)u^* \subseteq B$, and $u\in
N(B,M)$.  Then also $\mathcal{U}_{k+1}'' \subseteq
\mathcal{N}(B,M)''.$ The claim follows.  Then $M \subseteq
\mathcal{N}(B,M)'',$ and since the other containment holds
trivially, $M=\mathcal{N}(B,M)''.$ That is, $B$ is a Cartan masa
in $M$. Finally, by our construction of $M,$ condition (iv) will
imply that $\phi$ maps $M^n$ into $M,$ since $\phi$ is separately
normal in each variable.

     We proceed to construct the algebras $M_k, B_k.$  Put $B_0 =
M_0\cap A$.  Assume that $M_1,\ldots,M_k,$ and $B_1,\ldots,B_k$
have been constructed, satisfying conditions (i)-(iv), specified
above. By Dixmier's approximation theorem, \cite[Theorem
8.3.5]{KR}, each sequence $\{y_{k,r}\}_{r=1}^\infty$ is inside a
von Neumann algebra $Q_0 \subseteq N$ (generated by a countable
set of unitaries) such that $K_{Q_0}^n(y_{k,r}) \cap \mathbb{C}1$
is nonempty. To get condition (iii), observe that since
$\mathcal{N}(A)''=N,$ Lemma \ref{lemma:slippery} applies and there
is a countable set of unitaries $\mathcal{U}_{k+1}$ satisfying the
desired condition. Since $\phi$ is normal in each variable, by
\cite[Theorem 4.4]{cohomology survey}, $\phi$ is jointly
$\|\cdot\|_2$-norm continuous when restricted to bounded balls in
$(M_k)^n$.  Since $M_k$ is separable, it follows that
$\phi((M_k)^n)$ generates a separable von Neumann algebra $Q_1$.
Finally, since for all $r \in \mathbb{N},  \mathbb{E}_A(y_{k,r})
\in K_A^w(y_{k,r})\cap B_k,$ there exists a set of unitaries
$\{u_m\}_{m=1}^\infty$ generating a von Neumann algebra $Q_2
\subseteq A$ such that $\mathbb{E}_A(y_{k,r}) \in
K_{Q_2}^w(y_{k,r})$, for all r.  This will give condition (i).  We
complete the construction by letting $M_{k+1}$ be the von Neumann
algebra generated by $Q_0, Q_1, Q_2, \mathcal{U}_{k+1}, M_k,$ and
$\mathbb{E}_A(M_k)$.
\end{proof}

     We are now in a position to compute the cohomology groups
$H^n(N,N)$, where $N$ is a general type $\rm{II_1}$ factor with a
Cartan subalgebra.  Note that for any such $N$, associated to each
finite set $F \subseteq N$ is the separable von Neumann algebra
$M_F \subseteq N$ which it generates. Thus every such $\rm{II_1}$
factor $N$ satisfies the hypothesis of Proposition
$\ref{prop:cartansep}.$ This leads to our main theorem.

\begin{thm}\label{theorem:cartancohomology}
Let $N$ be a type $\rm{II_1}$ factor with a Cartan subalgebra $A.$
Then $H^n(N,N)=0$ for all $n \geq 1$.
\end{thm}
\begin{proof} Since the case $n=1$ is the Kadison-Sakai
result, we assume $n$ is at least 2.  We refer the reader to the
proof of the separable case in \cite{cartan sep}.  Let $\phi: N^n
\rightarrow N$ be a cocycle, which we may assume to be separately
normal. By Proposition $\ref{prop:cartansep}$, for each finite set
$F \subseteq N$, there exists a separable subfactor $N_F$ such
that $M_F \subseteq N_F \subseteq N$, $N_F \cap A$ is Cartan in
$N_F$, and $\phi$ maps $(N_F)^n$ into $N_F$. Denote the
restriction of $\phi$ to $N_F$ by $\phi_F$.  By the separable
case, $\phi_F$ is a coboundary, i.e., there exists an
$(n-1)$-linear and bounded map $\psi: (N_F)^{n-1} \rightarrow N_F$
such that $\phi_F=\partial \psi_F$.  Moreover, there is a uniform
bound on $\|\psi_F\|$; we confine this argument to the end of the
proof. Let $\mathbb{E}_F: N \rightarrow N_F$ be the conditional
expectation. Define $\theta_F:N^{n-1} \rightarrow N$ by
$\theta_F=\psi_F \circ (\mathbb{E}_F)^{n-1}.$ Order the finite
subsets of $M$ by inclusion.  Because $\|\psi_F\|$ is uniformly
bounded, so is $\|\theta_F\|$.  Now, for any $(n-1)$-tuple
$(x_1,\ldots,x_{n-1}) \in N^{n-1}$,
$\{\theta_F(x_1,\ldots,x_{n-1})\}_F$ is a bounded net as $F$
ranges over all finite sets containing $x_1,\ldots, x_{n-1}$.
This has an ultraweakly convergent subnet (which we also denote by
$\{\theta_F(x_1,\ldots,x_{n-1})\}_F$), by ultraweak compactness of
bounded subsets of $N$ (\cite[Chapter 7]{KR}). Define
$\theta:N^{n-1} \rightarrow N$ by
\[\theta(x_1,\ldots,x_{n-1}) = \lim_F \theta_F(x_1,\ldots,x_{n-1}). \]
Then $\theta$ is clearly $(n-1)$-linear, and bounded by the
uniform bound on $\|\theta_F\|$.  We claim $\phi=\partial \theta$.
Let $(x_1,\ldots,x_n) \in M^n$.  Then
\[\phi(x_1,\ldots,x_n)=\phi_F(x_1,\ldots,x_n)=\partial
\theta_F(x_1,\ldots,x_n),\] for all finite sets $F$ containing
$x_1,\ldots,x_n$.  Ordering these sets by inclusion, we will
obtain a subnet $\{\theta_F\}$ such that $\partial
\theta_F(x_1,\ldots ,x_n)$ converges weakly to $\partial
\theta(x_1,\ldots,x_n)$.  Passing to limits in the above equality
gives $\phi(x_1,\ldots,x_n)=\partial \theta(x_1,\ldots,x_n)$. This
proves the claim, and the result follows.

     It remains to be shown that there is a uniform bound on the
norms of the maps $\psi_F$, constructed above.  It suffices to
obtain an estimate for each of these maps in terms of $n$ and
$\|\phi_F\|$, since this last quantity is dominated by $\|\phi\|$
for all $F$.  We drop the index $F$, since it plays no further
role in the proof.

     Now $M$ will denote a separable $\rm{II_1}$ factor with Cartan
masa $A$.  By \cite[Theorem 2.2]{cartan sep}, there is a
hyperfinite factor $R$ such that \[A\subseteq R \subseteq M
\textrm{ and } R' \cap M=\mathbb{C}1.\] Let $\phi:M^n \rightarrow
M$ be a separately normal cocycle.  We will show that $\phi$ is
the image under $\partial$ of an $(n-1)$-linear map whose norm is
at most $K\|\phi\|$, where $K$ is a constant depending only on
$n$. By Lemma \ref{normal R modular} there exists a map $L_n$ such
that $\theta= \phi -
\partial L_n \phi$ is a separately normal $R$-module map, and $\|\theta\| \leq l(n) \|\phi\|$, where $l(n)$ is a
constant depending only on $n$.  Let $\mathcal{U}$ be a generating
set of unitaries for $M$.  Then the weak closure of
$C^*(\mathcal{U})$ is $M$. Let $\phi_C:C^*(\mathcal{U})^n
\rightarrow M$ denote the restriction of $\phi$. The proof of the
separable case gives $\theta= \partial \alpha$, where $\alpha:
C^*(\mathcal{U})^{n-1} \rightarrow M$ is $(n-1)$-linear and has
norm at most $\sqrt{2}\|\theta\|$.  Then
\[\phi_C= \theta + \partial \zeta = \partial (\alpha + \zeta),\]
where $\zeta$ denotes the restriction of $L_n \phi$ to
$C^*(\mathcal{U})^{n-1}$.  Write $\psi$ for $\alpha + \zeta$. Then
by what we have done above, there exists a constant $K(n)$ such
that \[\|\psi\| \leq K(n) \|\phi\|.\] We now wish to extend the
equality $\phi_C =\partial \psi$ to one involving maps defined on
$M$, while preserving this last norm estimate.  Applying the map
$V_n$ from Lemma \ref{extension2} to both sides of this last
equality, we get
\[V_n \phi_C=V_n (\partial \psi) = \partial (V_{n-1} \psi).\]
Because the separately normal map $\phi_C$ extends uniquely to
$\phi$ on $M$, this says $\phi = \partial (V_{n-1} \psi).$  The
map $V_{n-1}\psi$ is in $\mathcal{L}^{n-1}_w(M,M),$ and
$\|V_{n-1}\psi\| \leq K(n) \|\phi\|.$  This gives the required
norm estimate, completing the proof of the theorem.
\end{proof}

\end{document}